\documentclass[12pt]{article}
\usepackage{amssymb}
\usepackage{amsmath}
\usepackage{euscript}

\begin{document}

\centerline{\large\bf ERGODIC PROPERTIES OF DISCRETE}
 \centerline{\large\bf DYNAMICAL SYSTEMS }
\centerline{\large\bf AND ENVELOPING SEMIGROUPS}

\medskip
\centerline{\bf A.V. Romanov}
\medskip
\centerline{\it National Research University --- Higher School of Economics}
\centerline{\it e-mail: av.romanov@hse.ru}
\bigskip

{\small {\bf Abstract.} For a continuous semicascade on a metrizable compact
set $\Omega $,  we consider the weak$^{*}$ convergence of generalized operator
ergodic means in ${\rm End}\, \, C^{*} (\Omega )$. We discuss conditions on the
dynamical system under which (a) every ergodic net contains a convergent
subsequence; (b) all ergodic nets converge; (c) all ergodic sequences converge.
We study the relationships between the convergence of ergodic means and the
properties of transitivity of the proximality relation on $\Omega$, minimality
of supports of ergodic measures, and uniqueness of minimal sets in the closure
of trajectories of a semicascade.

These problems are solved in terms of three algebraic-topological
objects associated with the dynamical system: the Ellis enveloping
semigroup, the K\"{o}hler operator semigroup  $\Gamma $, and the
semigroup $G$ that is the weak$^{*} $ closure of the convex hull of
$\Gamma $ in ${\rm End}\, C^{*} (\Omega )$. The main results are
stated for ordinary semicascades (whose Ellis semigroup is
metrizable) and tame semicascades. For a dynamics, being ordinary
is equivalent to being ``nonchaotic'' in an appropriate sense. We
present a classification of compact dynamical systems in terms of
topological properties of the above-mentioned semigroups.

\noindent

\textit{Keywords}: ergodic property, nonchaotic dynamics, ordinary semicascade,
tame dynamical system, enveloping semigroup, Choquet representation.

\textit{2010 Mathematics Subject Classification.} Primary 37A30; Secondary
20M20.}

\bigskip

\centerline{ \bf Introduction}

\bigskip

The main topic of this paper is the weak$^{*} $ convergence of
generalized ergodic means along trajectories of discrete compact
dynamical systems. In the course of our study, ergodic properties
of the semicascade $(\Omega ,\varphi )$ generated by a continuous
transformation  $\varphi $ of a metrizable compact set $\Omega $
are related to characteristics of the Ellis enveloping semigroup
$E(\Omega ,\varphi )$, the K\"{o}hler operator enveloping semigroup
$\Gamma (\Omega ,\varphi )$, and yet another operator semigroup
$G(\Omega ,\varphi )$, which was considered in the author's paper
[1]. We also study relationships between properties of these three
semigroups. A considerable part of the results of the present paper
were announced in [2, 3].

The transformation $\varphi : \Omega \to \Omega $ induces linear
contractions $U$ and $V=U^{*} $ on the space $X=C(\Omega )$ of
continuous scalar functions on $\Omega $ and on the dual space
$X^{*} $ (consisting of Borel measures on $\Omega $), respectively.
Let ${\rm End}\, X$ and ${\rm End}\, X^{*} $ be the normed spaces
of bounded linear operators on $X$ and $X^{*} $, respectively. The
enveloping semigroup $E(\Omega ,\varphi )$ of the semicascade
$(\Omega ,\varphi )$ is defined [4] as the closure of the family
$\{ \varphi ^{n} , n\ge 0\} $ in the Hausdorff topology of
pointwise convergence on the space of all mappings $\Omega \to
\Omega $. The semigroup  $\Gamma (\Omega ,\varphi )$ is defined [5]
as the closure of $\Gamma _{0} =\{ \, V^{n} ,\, n\ge 0 \} $ in the
weak$^{*} $ operator topology ${\rm W}^{*} {\rm O}$ on the operator
space ${\rm End}\, X^{*} $. The operator semigroup $G(\Omega
,\varphi )$ is defined [1] as the ${\rm W}^{*} {\rm O}$-closed
convex hull of $\Gamma _{0} $. The semigroups $E(\Omega ,\varphi
)$, $\Gamma (\Omega ,\varphi )$, and $G(\Omega ,\varphi )$ are
compact in the respective topologies, and their topological
characteristics, as shown in [1], are closely related to ergodic
properties of the semicascade $(\Omega ,\varphi )$. It turns out
that such relationships manifest themselves especially strongly in
two classes of discrete dynamical systems. These are systems with
metrizable enveloping semigroup  $E(\Omega ,\varphi )$, dubbed
\textit{ordinary} in [2, 3], and \textit{tame} semicascades,
defined by the relation ${\rm card}\, E(\Omega ,\varphi ) \le
\mathfrak {c} $.

The main results of the paper pertain to the convergence of
generalized ergodic means for semicascades of the two
above-mentioned classes in the weak$^{*} $ operator topology ${\rm
W}^{*} {\rm O}$ on ${\rm End}\, X^{*} $. Let $G_{0}$ be the convex
hull of the set $\Gamma _{0} $ in ${\rm End}\, X^{*} $. We say that
a net of operators $T_{\alpha} \in G_{0} $ is ergodic (averaging)
if $T_{\alpha } (I-V)\stackrel{{\rm W}^{*} {\rm
O}}{\longrightarrow}0,\; \, I={\rm id}$. Then $T_{\alpha }
=R_{\alpha}^{*} $ for some $R_{\alpha } \in {\rm End}\, X$, and the
corresponding operator net $R_{\alpha }^{} $ consisting of the
elements of the convex hull of the operator set $\{ \, U^{n} ,\, \,
n\ge 0\, \} $ in End$\, X$ is said to be ergodic as well. The
sequences of Ces\`aro means
\[U_{n}   =  \frac{1}{n} \, (I+U+...+U^{n-1} ),\; \; \; \; \; \;
V_{n}  =  \frac{1}{n} \, (I+V+...+V^{n-1} ) \] are ergodic, since
$V_{n} (I-V) = n^{-1} (I-V^{n} )$,   $\left\| \, V^{n} \, \right\|
=1$, and $V_{n} =U_{n}^{*} $. The limit element of a convergent
ergodic operator net $T_{\alpha } \in G_{0} $ is always a bounded
projection on $X^{*} $. We point out that the  ${\rm W}^{*} {\rm
O}$-convergence of such nets is equivalent to the convergence, for
arbitrary continuous functions $x\in X$, of the means $R_{\alpha
} x\to z$, $z\in X^{**}, $ in the $w^{*} $-topology of the second dual
space $X^{**} = C^{**} (\Omega )$. For ergodic operator
\textit{sequences} $T_{n} \in G_{0} $, this convergence is
equivalent to the pointwise convergence of the functions $R_{n}
x\to z$ on $\Omega $, but in the general case it can be stronger
than the latter. Note that the idea of considering ${\rm W}^{*}
{\rm O}$-convergence of ergodic means essentially goes back to the
classical paper [6] by Krylov and Bogolyubov.

Issues related to the ${\rm W}^{*} {\rm O}$-convergence of
arbitrary ergodic operator nets or their certain subsequences are
considered here in connection with purely dynamic properties of
ordinary and tame dynamical systems: transitivity of the
proximality relation, uniqueness of a minimum set in the closure of
an arbitrary trajectory, and minimality of supports of ergodic
measures. The techniques used here include the properties of the
enveloping semigroups $E(\Omega ,\varphi )$, $\Gamma (\Omega
,\varphi )$, and $G(\Omega ,\varphi )$ associated with the
semicascade $(\Omega ,\varphi )$. In the ordinary case, the key
point is the metrizability of the semigroup \textbf{$G(\Omega
,\varphi )$}, which permits one to pass from the convergence of
ergodic nets to that of ergodic sequences. By way of example, let
us present some of the results obtained below for the ${\rm W}^{*}
{\rm O}$-convergence of ergodic means.

(i) For an ordinary semicascade \textbf{$(\Omega ,\varphi )$},
every ergodic operator net $T_{\alpha } $ contains a convergent
subsequence $T_{\alpha (k)} $. Moreover, all ergodic operator nets
$T_{\alpha} $ converge if and only if, for each $\omega \in \Omega
$, the closure of the trajectory $\{ \varphi ^{n} \omega ,\; n\ge
0\} $ contains a unique minimal set. The same is true for sequences
instead of nets.

(ii) If, for a tame semicascade \textbf{$(\Omega ,\varphi )$}, the
closure of each trajectory $\{ \, \varphi ^{n} \omega ,\, \, n\ge
0\, \} ,\; \, \omega \in \Omega {\kern 1pt}$ contains a unique
minimal set, then all ergodic operator sequences $T_{n} $ converge,
and the support of each $\varphi $-ergodic measure is a minimal
set.

   As was shown in [1], for an \textit{arbitrary}
semicascade $(\Omega ,\varphi )$, the convergence of all ergodic
operator nets is equivalent to the condition ${\rm card}\, L=1$,
where $L$ is the kernel (the intersection of all two-sided ideals)
of the semigroup $G(\Omega ,\varphi )$; it consists of the
operators $Q\in G(\Omega ,\varphi )$ with the property $VQ=Q$. The
convex weakly$^{*} $ compact set $L\subset {\rm End}\, X^{*} $ of
projections was originally considered by Lloyd [7]. There are
reasons to assume that it is the algebraic-geometric
characteristics of this object that are responsible for the ${\rm
W}^{*} {\rm O}$-ergodic properties of the semicascade $(\Omega
,\varphi )$. Let us also mention one more significant result
related to purely dynamic properties of tame semicascades: for
these semicascades, the closure of the union of all minimal sets
coincides with the minimal attraction center, i.e., with the
closure of the union of supports of all $\varphi $-ergodic Borel
measures on $\Omega $. In this connection, note the strict
ergodicity of minimal tame systems, which has recently been
established in [8, 9].

Numerous properties of Ellis semigroups of ordinary and tame
systems (not necessarily discrete) can be found in Glasner's survey
[10]. As follows from the results in [10--12], a semicascade
$(\Omega ,\varphi )$ is ordinary if and only if the dynamics is
nonchaotic in the sense that each closed semi-invariant ($\varphi
\, \Theta \subset \Theta $) set $\Theta \subset \Omega \; $
contains a trajectory that is Lyapunov stable with respect to the
restriction $(\Theta ,\varphi )$. It is also well known [10] that
all weakly almost periodic dynamical systems prove to be ordinary.

The present paper suggests an alternative definition for $(\Omega
,\varphi )$ bering ordinary, which is the requirement that the
semigroup $G(\Omega ,\varphi )$ is metrizable. Moreover, we show
that, in the ordinary case, the convex compact subset $G(\Omega
,\varphi )$ of the linear space ${\rm End}\, X^{*} $ is a Choquet
simplex, and each operator $T\in G(\Omega ,\varphi )$ is uniquely
determined by the values it takes on the Dirac measures $\delta
(\omega )$, $\omega \in \Omega $.

The class of tame dynamical systems, originally introduced (under a
different name) by K\"{o}hler [5], was comprehensively studied in
[10--13]. In the present paper, we speak of semicascades $(\Omega
,\varphi )$ satisfying any of the following three equivalent
conditions: (1) the semigroup $E(\Omega ,\varphi )$ is at most of
the cardinality of continuum; (2) the compact set $E(\Omega
,\varphi )$ is a Fr\'echet--Urysohn space; (3)  all transformations
$p\in E(\Omega ,\varphi )$ belong to the first Baire class, i.e.,
are pointwise limits of sequences of continuous mappings $p_{n}
:{\kern 1pt} {\kern 1pt}  \Omega \to \Omega $. We see that the
class of tame dynamical systems includes the ordinary semicascades.
Here we establish one more criterion for tame dynamics: the
enveloping semigroup $E(\Omega ,\varphi )$ consists of Borel
transformations.

The correspondence  $V\, \, \to \, \, \varphi $ generates a
continuous algebraic semigroup epimorphism  $\Gamma (\Omega
,\varphi )\stackrel{\pi}{\longrightarrow}E(\Omega ,\varphi )$ [5].
Moreover, $P\delta (\omega )=\delta (p\, \omega )$  for $P\in
\Gamma (\Omega ,\varphi )$,  $p=\pi P$,  $\omega \in \Omega$, and
the Dirac measures $\delta (\cdot )$. Let $A(\Omega )$ be the set
of Borel probability measures on the compact set $\Omega $. It is
well known [5] that, for a tame system $(\Omega ,\varphi )$, the
function $\pi $ is injective and the semigroup $\Gamma (\Omega
,\varphi )$ can be identified with $E(\Omega ,\varphi )$. It turns
out (Corollary 2.6) that in this case all operators $P\in \Gamma
(\Omega ,\varphi )$ act on the measures $\mu \in A(\Omega )$ by the
``natural'' rule $(P\mu )(e)=\mu (p^{-1} e)$ for $p=\pi P$ and
Borel sets $e\subset \Omega $  (i.e., by analogy with the operators
$V^{n} , n\ge 0$, generating the semigroup $\Gamma (\Omega ,\varphi
)$). It is also proved that, for the case of a tame semicascade
$(\Omega ,\varphi )$, the set of extreme points of the convex
compact subset $G(\Omega ,\varphi )$ of the operator space ${\rm
End}\, X^{*} $ coincides with the set $\Gamma (\Omega ,\varphi )$.

We point out that although Ellis' monograph [4] deals with
invertible dynamical systems, the results in [4] cited in what
follows are easily seen to remain valid in the noninvertible case.

The subsequent exposition is organized as follows. Section 1
provides general information about the enveloping semigroups
$E(\Omega ,\varphi )$, $\Gamma (\Omega ,\varphi )$, and $G(\Omega
,\varphi )$ and describes the relationships between these
semigroups. Section 2 contains a broad hierarchical compendium (in
a sense, a classification) of earlier-known and some new
dependences between various properties of arbitrary semicascades,
the corresponding statements being mostly given in terms of the
above-mentioned semigroups. The main results pertaining to
ergodicity properties of ordinary and tame semicascades are
presented in Sec.~3.

\bigskip

 \centerline{ \bf1. Preliminaries}

\bigskip

Let $\varphi $ be a continuous (not necessarily invertible)
transformation of a metrizable compact set $\Omega $, and let
$(\Omega ,\varphi )$ be the corresponding discrete dynamical system
(semicascade) on $\Omega $. We define the shift operator
$U=U_{\varphi} $ on the space $X=C(\Omega )$ of continuous scalar
functions by the formula $(Ux)(\omega )=x(\varphi \omega )$, $x\in
C(\Omega )$. By $A(\Omega )$ and $K(\Omega )$ we denote the sets of
Borel probability measures and Dirac measures on $\Omega $,
respectively. Recall that every Borel measure on a metrizable
compact set is regular. Also, let $\Lambda (\Omega )$ be the class
of ergodic measures of the semicascade $(\Omega ,\varphi )$, i.e.,
$\varphi $-invariant measures $\lambda \in A(\Omega )$ such that
$\lambda (e)$ is either~$0$ or~$1$ for every Borel set $e\subset
\Omega $ with $\varphi ^{-1} e=e$. The set $A(\Omega )$ is convex,
compact, and metrizable in the $w^{*} $-topology of the dual space
$X^{*} =C^{*} (\Omega )$, and $K(\Omega )$ is a closed subset of
$A(\Omega )$.

By $ \rm Im\, (\cdot )$ and $F(\cdot )$ we denote the range and the
subspace of fixed vectors, respectively, of a linear operator.
Next, ${\rm End}\, X$  and ${\rm End}\, X^{*}$ are the Banach
algebras of bounded linear operators on $X$ and on the first dual
space $X^{*} $, respectively. Let $V=U^{*} $ be the adjoint
operator on $X^{*} $; then $\left\| \, V^{n} \, \right\| =1$ for
$n\ge 0$. We equip the operator algebra ${\rm End}\, X^{*} $ with
the locally convex weak$^{*} $ operator topology (${\rm W}^{*} {\rm
O}$-topology) and single out the multiplicative semigroups $\Gamma
_{0} =\{V^{n} ,\, n\ge 0 \}$ and $G_{0} = {\rm co}\, \Gamma _{0}$
in this algebra. Set $G=\overline{\rm co} \, \Gamma _{0}$ and
$\Gamma =\overline{\Gamma}_{0}$, where $\overline{(\,\cdot \,)}$ is
the operation of ${\rm W}^{*} {\rm O}$-closure of sets in the
normed space ${\rm End}\, X^{*} $; then $\Gamma \subset G$ and
$\left\| \, T\, \right\| =1$ for $T\in G$. Next, $T: A(\Omega )\to
A(\Omega )$ for the operators $T\in G$. Norm-bounded subsets of the
space ${\rm End}\, X^{*} $ are relatively compact in the ${\rm
W}^{*} {\rm O}$-topology, and hence the sets $\Gamma =\Gamma
(\Omega ,\varphi )$ and $G=G(\Omega ,\varphi )$ are Hausdorff
(generally, nonmetrizable) separable compact sets. The
multiplication in the algebra ${\rm End}\, X^{*} $ is only
continuous in the left argument; as a consequence, the semigroup
$G$ and the subsemigroup $\Gamma $ are usually noncommutative,
although their elements commute with the operators $V^{n} , \, n\ge
0$. Note that ${\rm Im}\, T \supset F(T)\supset F(V)$ for all $T
\in  G$. As was shown in [1], the kernel (the intersection of all
two-sided ideals) of $G(\Omega ,\varphi )$ coincides with the
closed set
 \[L \, =\,  \{ \, Q\in G:\, VQ=Q\, \} .
                       \eqno           (1.1)\]
Furthermore, $TQ=Q$ for all $T \in  G$ and  $Q \in L$, and the
convex compact set $L$ consists of unit-norm projections. One can
give an alternative definition, $L \, =\,  \{  T \in  G:\, {\rm
Im}\, T= F(V) \} $.

The notion of enveloping semigroup of a dynamical system plays an
important role in what follows. According to [4], the enveloping
semigroup $E(\Omega ,\varphi )$ of a semicascade $(\Omega ,\varphi
)$ is the closure of the generating family $E_{0} (\Omega ,\varphi
)=\{ \varphi ^{n} ,\, n\ge 0\} $ in the Hausdorff topology $\tau $
of pointwise convergence on the space $\Omega ^{\Omega } $ of all
mappings $p:\Omega \to \Omega $. The semigroup $E(\Omega ,\varphi
)$ is compact but nonmetrizable and noncommutative in general, and
its center contains the set $E_{0} (\Omega ,\varphi )$. Every
enveloping semigroup is separable and has nonempty kernel.

Since $(V\mu )(e)=\mu (\varphi ^{-1} e)$ for $\mu \in A(\Omega )$
and Borel sets $e\subset \Omega $, it follows that the operator $V$
generates the semicascade $(A,V)$ on the compact set $A=A(\Omega
)$. According to [5], the semigroup $\Gamma (\Omega ,\varphi )$ can
be identified with the enveloping semigroup $E(A,V)$ of $(A,V)$,
which is the closure of the family $\{ \, V^{n} , \, n\ge 0 \} $ of
continuous transformations of the compact set $A(\Omega )$ in the
Cartesian product topology on $A^{A} $. On the other hand, it is
convenient to interpret the semigroup $G(\Omega ,\varphi )$ as the
enveloping semigroup of the action $W\times A\,
\stackrel{V}{\longrightarrow} A$ on $A=A(\Omega )$ of the abelian
semigroup $W$ of nonnegative finite numerical sequences with unit
sum and with convolution as multiplication. The semigroup of comvex
combinations of the monomials $\{ \, t^{n} ,\, n\ge 0 \} $ with the
usual multiplication may serve as an alternative model for~$W$.
Thus, $G(\Omega ,\varphi )\simeq E(A,W)$.

Let $\Sigma _{b} $ be the $\sigma$-algebra of Borel subsets of
$\Omega $, and let $X_{b} $ be the family of scalar bounded Borel
functions on $\Omega $. Clearly, $X\subset X_{b} \subset X^{**} $,
and $X_{b} $ is a norm-closed subspace of the second dual $X^{**}$.
Moreover $\left\| \, x\, \right\| =\sup \, \left|\; x(\omega )\,
\right|$, where the supremum is over $\omega \in \Omega $, for
$x\in X_{b} $. Note that Borel functions are universally
measurable, i.e., measurable with respect to the Lebesgue extension
of an arbitrary measure  $\mu \in A(\Omega )$. We say that a
mapping $p: \Omega \to \Omega $ is a \textit{Borel} mapping if
$p{}^{-1} e\in \Sigma _{b} $ for all closed sets $e\subset \Omega
$. Finally, let $\Pi _{b} $ be the set of all Borel transformations
of $\Omega $. Note that actually $p{}^{-1} e\in \Sigma _{b} $ for
$p\in \Pi _{b} $ and $e\in \Sigma _{b} $.

Next,  $V\delta (\omega )=\delta (\varphi \, \omega )$ for $\omega
\in \Omega $, and the set $K=K(\Omega )$ of Dirac measures is
closed in the topological space $A(\Omega )$; hence $VK\subset K$
and $PK\subset K$ for all operators $P\in \Gamma (\Omega ,\varphi
)$. The correspondence $V \to \varphi $ generates a continuous
algebraic homomorphism $\pi :\, P\to p$ of the semigroup $\Gamma
(\Omega ,\varphi )$ into the semigroup $E(\Omega ,\varphi )$, whose
action is clearly defined by the relation
  \[P\delta (\omega )=\delta (p\, \omega )\, ,\quad
  \omega \in \Omega \, .   \eqno    (1.2)\]
Since the set $\Gamma (\Omega ,\varphi )$ is compact, $\pi V^{n}
=\varphi ^{n}, $  and the generating transformation family $E_{0}
(\Omega ,\varphi )$ is $\tau $-dense in $E(\Omega ,\varphi )$, we
can speak of the \textit{canonical epimorphism} $\pi :\; \Gamma
(\Omega ,\varphi )\to E(\Omega ,\varphi )$.

For a measure $\mu  \in  A(\Omega )$ and a set $e \in  \Sigma
_{b}$,  we have $(V^{n} \mu )(e) = \mu (\varphi ^{-n} e)$, $n \ge
0$. Moreover, $(V^{n} )^{*} : \, x(\omega )\to x(\varphi ^{n}
\omega )$ for continuous functions $x\in X$  and for  $\omega \in
\Omega $. It is useful to find out in what cases arbitrary
operators $P\in \Gamma (\Omega ,\varphi )$ and their duals $P^{*}
\in {\rm End}\, X^{**} $ act in a similar way.

DEFINITION 1.1.  We say that an operator $P\in \Gamma
(\Omega,\varphi )$ is \textit{regular} if  $p\in \Pi _{b} $ for
$p=\pi P$ and
 \[(P\mu )(e) = \mu (p^{-1} e)    \eqno          (1.3)\]
for  $\mu \in A(\Omega )$ and $e\in \Sigma _{b} $.

LEMMA 1.2.\textit{  If  $P\in \Gamma (\Omega ,\varphi )$ and $p=\pi
P$}, \textit{then the following conditions are equivalent}:

(a)\textit{   $P^{*} :\,  X\to X_{b} $}.

(b)\textit{   $(P^{*} x)(\omega )=x(p\omega )$  for  $x\in X$ and}
$\omega \in \Omega $.

(c)\textit{ The operator $P$ is regular}.

PROOF.  Clearly,  (b)  $\Rightarrow $  (a).  If condition (c)
holds, then $(x,P\mu )=(x(p\, \omega ),\, \mu )$ for any $x \in X$
and $\mu \in A(\Omega )$. On the other hand, one always has $(P^{*}
x,\mu )=(x,P\mu )$, and hence $(P^{*} x)(\omega ) \, =x(p\, \omega
)$ and (c) $\Rightarrow $ (b).  Now let us establish the
implication (a) $\Rightarrow $ (c). Let $x \in  X$, $y=P^{*} x$ and
$\omega \in \Omega $; then $y\in X_{b} $ and $(y,\delta (\omega
))=x(p\, \omega )$. The elements $y\in X_{b} $, being functionals
in $X^{**} $, are completely determined by their values on the
Dirac measures, and hence $(P^{*} x)(\omega )=x(p\omega )$ for all
$\omega \in \Omega $. Thus, $x(p\omega )\in X_{b} $ for each
continuous function $x(\omega )$. The indicator function $x_{e}
(\omega )$ of an arbitrary closed set $e\subset \Omega $ is the
pointwise limit of some sequence of continuous functions $x_{n}
(\omega )$; consequently, $x_{n} (p\, \omega )\to x_{e} (p\, \omega
)$ for all $\omega \in \Omega $, and hence the function $x_{e}
(p\omega )$ is a Borel function. It follows that $p\in \Pi _{b} $.
Next,
   \[(P^{*} x,\mu )=\int _{\Omega }x(p\, \omega )\, \mu (d\omega )
    =\int _{\Omega }x(\omega )\, \mu _{p} (d\omega )=(x,P\mu ) \]
for any $\mu \in A(\Omega )$, where the Borel probability measure
$\mu _{p} $ is defined on the sets  $e\in \Sigma _{b} $ by the
formula $\mu _{p} (e) = \mu (p^{-1} e)$. It follows that $P\mu =\mu
_{p} $, and the proof of the lemma is complete.

We refer to arbitrary finite linear combinations of Dirac measures
as discrete measures. Since Borel functions are universally
measurable, we can identify the dual space $X^{*} =C^{*} (\Omega )$
with a strongly closed subspace of $X_{b}^{*} $. Let
\textbf{$\sigma $} be the topology on $X^{*} $ corresponding to the
convergence of linear functionals on the functions $x\in X_{b} $.
Clearly, \textbf{$\sigma \ge w^{*} $} in the sense of comparison of
topologies. Recall that an operator $P\in \Gamma (\Omega ,\varphi
)$ acting on $X^{*} $ is continuous in the \textbf{$w^{*}
$}-topology only if $P=R^{*} $ for some $R\in {\rm End}\, X$. At
the same time, the regularity of an operator $P\in \Gamma (\Omega
,\varphi )$ implies its \textbf{$(\sigma ,w^{*} )$}-continuity,
i.e., continuity with the \textbf{$\sigma $}-topology on the domain
and the \textbf{$w^{*} $}-topology on the range. This readily
follows from the identity \textbf{$(x,Py)=(P^{*} x,y)$} for
\textbf{$x\in X$} and \textbf{$y\in X^{*} $ } and from the fact
that \textbf{$P^{*} : X\to X_{b} $} by Lemma~1.2.

LEMMA 1.3.  \textit{The linear manifold of discrete measures is
$\sigma $-dense in} $X^{*} $.

PROOF.  The Dirac measures form a total set in the dual space
$X_{b}^{*} $. Consequently, according to general facts about weak
topologies [14, Chap.~5, Sec.~12], the linear manifold generated by
these measures is $\sigma $-dense in $X_{b}^{*} $. Thus, the
discrete measures are also $\sigma $-dense in the space $X^{*} $ of
functionals. The proof of the lemma is complete.

Note that every Borel measure on a metrizable compact set is a
sequential \textbf{$w^{*} $-}limit of discrete measures. (This
follows from the Riemann integrability of continuous functions.)

Set  $\Phi (P)=\pi ^{-1} \pi P$  for  $P\in \Gamma (\Omega ,\varphi
)$.  The operator classes $\Phi (P)$ form a partition of $\Gamma
(\Omega ,\varphi )$, which is generated by the equivalence relation
``$P_{1} \sim P_{2} $ whenever $\pi P_{1} =\pi P_{2} $.''

LEMMA 1.4. \textit{If $P\in \Gamma (\Omega ,\varphi )$},
\textit{then the set $\Phi (P)$ contains at most one regular
element}.

PROOF.  If $p=\pi P$ for some regular operator $P\in \Gamma (\Omega
,\varphi )$, then $p\in \Pi _{b} $. By Lemma 1.2, $(P^{*} x)(\omega
)=x(p\omega )$ and hence also $(x,P\delta (\omega ))=(P^{*}
x,\delta (\omega ))=x(p\omega )$ for any  $x\in X$ and $\omega \in
\Omega $. By Lemma 1.3, the discrete measures are $\sigma $-dense
in the dual space $X^{*} $, and so it follows from the
\textbf{$(\sigma ,w^{*} )$-}continuity of the regular operator
$P\in {\rm End}\, X^{*} $ that the restriction of its adjoint
$P^{*} \in {\rm End}\, X^{**} $ to the subspace $X$, which is
$w^{*} $-dense in $X^{**} $, is determined by a Borel
transformation $p: \Omega \to \Omega $. The operator $P^{*}: \,
X^{**} \to X^{**} $ is continuous in the \textbf{$w^{*} $}-topology
on $X^{**} $ and hence is completely determined by its values on
$X$ and thus by $p$. The proof of the lemma is complete.

Let ${\rm sq}\overline{(\,\cdot \,)}$ be the set of all possible
$\tau $-limits of sequences of elements in an arbitrary set of
transformations of $\Omega $. Next, let  $E_{b} (\Omega ,\varphi )$
be the least sequentially $\tau $-closed subset of all
transformations $p: \Omega \to \Omega $ containing the generating
family $E_{0} (\Omega ,\varphi )$. In more detail, let us construct
an increasing transfinite sequence of transformation classes
$E_{\nu } =E_{\nu } (\Omega ,\varphi )$ by setting $E_{\nu }
=\bigcup _{\eta <\nu}E_{\eta }  $ or  $E_{\nu } ={\rm sq}(E_{\nu
-1} )$ depending on whether  $\nu $ is a limit ordinal number or
not. Then the union $E_{b} (\Omega ,\varphi )$ of classes $E_{\nu}
$ over all ordinal numbers $\nu $ can be obtained from the
generating transformation family $E_{0} (\Omega ,\varphi )$ by a
transfinite sequence of sequential passages to the $\tau $-limit.
The set $\Pi _{b} $ of Borel transformations is stable under such
passages, and $E_{0} (\Omega ,\varphi )\subset \Pi _{b} $, so that
$E_{b} (\Omega ,\varphi )\subset E(\Omega ,\varphi )\bigcap \Pi
_{b} $; moreover, the inclusion can be proper.

The following claim (which is of interest in itself) is used in the
derivation of Theorem 2.5 and Corollary 2.6 in what follows.

PROPOSITION 1.5.  \textit{If the function $\Gamma (\Omega ,\varphi
)\stackrel{\pi}{\longrightarrow}E(\Omega ,\varphi )$ defined in
}(1.2)\textit{  is injective and $E(\Omega ,\varphi )=E_{b} (\Omega
,\varphi )$},\textit{  then all operators $P\in \Gamma (\Omega
,\varphi )$ are regular}.

The proof will be given below.

Let ${\rm End}\, (X,X^{**} )$ be the Banach space of bounded linear
operators from $X$ to $X^{**} $. Suppose that operators $R\in {\rm
End}\, (X,X^{**} )$ and $T\in {\rm End}\, X^{*} $ are related by
the identity $(Rx,y)=(x,Ty)$ for $x\in X$ and $y\in X^{*} $. Then
$R=\left. T^{*} \right|_{X} $ and $T=\left. R^{*} \right|_{X^{*} }
$, and moreover, the correspondence
$R\stackrel{\theta}{\longrightarrow}T$ defines a natural linear
isometry between the spaces ${\rm End}\, (X,X^{**} )$ and ${\rm
End}\, X^{*} $. To each transformation $p\in \Pi _{b} $, we assign
a linear operator $R\in {\rm End}\, (X,X^{**} )$ by setting
$(Rx)(\omega )=x(p \, \omega )$ for  $x\in X$ and $\omega \in
\Omega $. Then $Rx\in X_{b} $ and $P\in {\rm End}\, X^{*} $ if we
set $P=\theta R=\psi p$, where the function $\psi :\,  \Pi _{b} \to
{\rm End}\, X^{*} $, $p\stackrel{\psi}{\longrightarrow}P$, is
generated by relation (1.3).

LEMMA 1.6.  \textit{If  $p\in E_{b} (\Omega ,\varphi )$ and $P=\psi
p$},\textit{  then} $P\in \Gamma (\Omega ,\varphi )$.\textbf{}

PROOF.  Clearly, $P\in \Gamma (\Omega ,\varphi )$  for $p\in E_{0}
(\Omega ,\varphi )$. Every element $p\in {\rm sq}(E_{\nu } )$ is
the $\tau $-limit of a sequence of elements $p_{k} \in E_{\nu
} $, and hence  $x(p_{k} \omega )\to x(p\omega )$ for all
$x\in X$,  $\omega \in \Omega $.  If we set $(R_{k} x)(\omega
)=x(p_{k} \omega )$ and $(Rx)(\omega )=x(p\omega )$, then $R_{k}
,\, \, R\in {\rm End}\, (X,X^{**} )$, and by Lebesgue's dominated
convergence theorem we have $(R_{k} x,\mu )\to (Rx,\mu )$ for $\mu
\in A(\Omega )$. It follows that $P_{k} \to P$ in the ${\rm W}^{*}
{\rm O}$-topology of the space ${\rm End}\, X^{*} $, where $P_{k}
=\theta R_{k} $ and $P=\theta R$. Thus, if $P_{k} \in \Gamma
(\Omega ,\varphi )$, then $P\in \Gamma (\Omega ,\varphi )$ as well.
An application of transfinite induction over ordinal numbers $\nu $
to the transformation classes $E_{\nu } (\Omega ,\varphi )$
completes the proof.

We see that $\pi \psi ={\rm id}$ on the set $E_{b} (\Omega ,\varphi
)$.

LEMMA 1.7.  \textit{Let  $P\in \Gamma (\Omega ,\varphi )$ and $\pi
P\in E_{b} (\Omega ,\varphi )$}. \textit{  Then the operator set
$\Phi (P)$ contains exactly one regular element}.

PROOF.  If $P_{0} =\psi \pi P$,  then Lemma 1.6 guarantees the
inclusion $P_{0} \in \Gamma (\Omega ,\varphi )$.  Moreover, $P_{0}
\in \Phi (P)$  and $P_{0}^{*} :\; X\to X_{b} $, and hence the
operator $P_{0} $ is regular by Lemma 1.2. The uniqueness of such
operators in the class $\Phi (P)$ is ensured by Lemma 1.4, and the
proof is complete.

Now Proposition 1.5 directly follows from Lemma 1.7, because the
proposition assumes that $\Phi (P)=\{ P\} $.

The continuous action $P\to VP$ generates the semicascade $(\Gamma
,V)$ on the ${\rm W}^{*} {\rm O}$-compact set  $\Gamma \subset {\rm
End}\, X^{*} $, where $\Gamma =\Gamma (\Omega ,\varphi )$. We adopt
the following notation: $A(\Gamma )$ is the compact (in the $w^{*}
$-topology on $C^{*} (\Gamma )$) convex set of regular Borel
probability measures on $\Gamma $; ${\rm Ai}(\Gamma )$ is the
closed convex subset  of $V$-invariant measures in $A(\Gamma )$;
$\Lambda (\Gamma )$ is the subset of $V$-ergodic measures in ${\rm
Ai}(\Gamma )$, i.e., by one of the equivalent definitions in [15,
Chap.~10], the set ${\rm ex}\, {\rm Ai}(\Gamma )$ of extreme
points.

Next, since $G=\overline{\rm co} \, \Gamma _{0}$ and $\Gamma
=\overline{\Gamma}_{0}$, where $\Gamma _{0} =\{  V^{n} ,\, n\ge 0
\} $, it follows according to [15, Chap.~1] that ${\rm ex}\,
G\subset \Gamma $, and since $\Gamma $ is closed, it follows  [15,
Chap. 4] that to each operator $T\in G$ there corresponds a measure
$\lambda \in A(\Gamma )$ realizing the Choquet representation
  \[T=\int _{\Gamma }P\, \lambda (dP) .
                               \eqno                (1.4)\]
The integration here can be understood in the sense that
 \[(x,T\mu )= \int _{\Gamma }(x,P\mu )\, \lambda (dP)
        \eqno       (1.5)\]
for all $x\in X$ and $\mu \in A(\Omega )$. We say that $\lambda $
is a \textit{representing} measure (for the operator $T$) and
sometimes write $T=T_{\lambda} $. Note that the representing
measure is not necessarily unique. The mapping $\gamma :\lambda \to
T_{\lambda} $ of $A(\Gamma )$ onto $G(\Omega ,\varphi )$ (the
Choquet epimorphism) is continuous with the  $w^{*} $-topology on
the domain and the ${\rm W}^{*} {\rm O}$-topology on the range.

It turns out that the correspondence $\lambda \to T_{\lambda } $
relates invariant measures on the compact set $\Gamma (\Omega
,\varphi )$ to elements of the kernel $L$ of the semigroup
$G(\Omega ,\varphi )$.

LEMMA 1.8.  \textit{The function  $\gamma $ maps ${\rm Ai}(\Gamma
)$ onto the set}  $L$.

PROOF.  Let  $\lambda \in \; {\rm Ai}(\Gamma )$;  then the relation
   \[\int _{\Gamma }g(VP)\, \lambda (dP)
   = \int _{\Gamma }g(P)\, \lambda (dP)          \eqno    (1.6)\]
holds for all functions $g\in C(\Gamma )$, including functions of
the form $g(P)=(x,P\mu )$ with given $x\in X$ and $\mu \in A(\Omega
)$. We find from (1.5) and  (1.6) that
    \[\int _{\Gamma }(x,VP\mu )\, \lambda (dP) = \int _{\Gamma }(Ux,P\mu )\, \lambda
     (dP)\, = \,  (Ux,T\mu )=(x,VT\mu )
     =(x,T\mu ),  \]
and hence $VT=T$ and $T\in L$. Next, the continuous transformation
$P\to VP$ of the compact set $\Gamma (\Omega ,\varphi )$ induces
the (also continuous) self-mapping $\upsilon $ of the convex
compact set $A(\Gamma )$ of measures. Moreover, for each operator
$Q\in L$ there exists a measure $\lambda \in A(\Gamma )$ such that
     \[VQ=\int _{\Gamma }P\, \lambda _{\, 1} (dP) ,
     \; \; \; \; \;   \lambda _{\, 1} =\upsilon \lambda .  \]
Since  $VQ=Q$, we see that this mapping preserves the convex
compact sets $\gamma ^{-1} Q\subset A(\Gamma )$, and its fixed
points $\lambda \in \gamma ^{-1} Q$ prove to be invariant measures
in ${\rm Ai}(\Gamma )$. The proof of the lemma is complete.

Thus, for each operator $Q\in L$ there exists at least one
representation (1.4) with $V$-invariant measure $\lambda $ on the
semigroup $\Gamma (\Omega ,\varphi )$. Now let us establish a
relationship between ergodic measures $\lambda \in \Lambda (\Gamma
)$ and extreme points of the kernel $L$ of the semigroup $G(\Omega
,\varphi )$.

LEMMA 1.9.  \textit{For each operator  $Q\in {\rm ex}\, L$, the
full preimage $\gamma ^{-1} Q$ contains an ergodic measure $\lambda
\in \Lambda (\Gamma )$}. \textit{If the Choquet epimorphism
$\gamma:{\rm Ai}(\Gamma )\to L$ is injective},\textit{ then}
$\gamma : \Lambda (\Gamma )\to {\rm ex}\, L$.

PROOF. Let $Q\in {\rm ex}\, L$; then $\gamma ^{-1} Q$ is an extreme
subset of ${\rm Ai}(\Gamma )$, and hence [14, Chap.~5, Sec.~8] the
extreme points of $\gamma ^{-1} Q$ are extreme points of ${\rm
Ai}(\Gamma )$, and ${\rm ex}\, \gamma ^{-1} Q \subset \gamma ^{-1}
Q \bigcap \Lambda (\Gamma )$. The opposite inclusion is always
true, and it follows that every element of $Q\in {\rm ex}\, L$ has
at least one ergodic representing measure. Now let $\lambda \in
\Lambda (\Gamma )$ and $Q=\gamma \lambda $. If $2Q=Q_{1} +Q_{2} $
with $Q_{i} \in L$, then, by Lemma 1.8, $Q_{i} =\gamma \lambda _{i}
$ for some measures $\lambda _{1}, \lambda _{2} \in {\rm Ai}(\Gamma
)$. The assumption that the restriction of $\gamma $ to ${\rm
Ai}(\Gamma )$ is injective ensures that $2\lambda =\lambda _{1}
+\lambda _{2} $, $\lambda =\lambda _{i} $, and $Q=Q_{i} $, $i=1,2$.
The proof of the lemma is complete.

It follows from formula (1.4) with regard for (1.2) that
  \[T_{\lambda } \delta (\omega )
  =\int _{\Gamma }\delta (p \, \omega )\, \lambda (dP)
  \; \; \; \; \; \; \; (\omega \in \Omega ),  \]
where  $\lambda \in A(\Gamma )$ and $T_{\lambda } \in G$. The
extreme points of the convex set $A(\Omega )\subset X^{*} $ exactly
coincide with the set $K(\Omega )$ of Dirac measures $\delta
(\omega )$, $\omega \in \Omega $, and hence the Choquet
representation for a measure $\mu \in A(\Omega )$ can be written as
   \[\mu = \int _{\Omega }\delta (\omega )\, \mu (d\omega ) .  \]
In particular, this is true for $\mu =T_{\lambda } \delta (\xi )$.
If $e\in \Sigma _{b} $ and $h=\{P\in \Gamma (\Omega ,\varphi ):
p\,\xi \in e\, \} $, where $p=\pi P$, then $\mu (e)=\lambda (h)$.
For given $\xi \in \Omega $, the corresponding $e\to h$ carries the
operations of union, intersection and complement of sets from
$\Omega $ to $\Gamma (\Omega ,\varphi )$ and preserves the classes
of closed sets and hence of Borel sets. Moreover,
  \[\{ P\in \Gamma :  \, VP\in h \} =\{  P\in \Gamma :\,
  \varphi p\,\xi \in e \} =\{  P\in \Gamma :\,  p\,\xi \in \varphi
 ^{-1} e \} ,  \]
whence it follows that the correspondence $\lambda \to T_{\lambda }
\delta (\xi )$ preserves the invariance of measures; i.e., ${\rm
Ai}(\Gamma )\to {\rm Ai}(\Omega )$. Finally, to a $\varphi
$-invariant set $e\subset \Omega $ there corresponds a
$V$-invariant set $h\subset \Gamma $; consequently, ergodic
measures $\lambda \in \Lambda (\Gamma )$ are taken to ergodic
measures $\mu \in \Lambda (\Omega )$. Thus, we arrive at

LEMMA 1.10.  \textit{If  $\omega \in \Omega $ and $\lambda \in
\Lambda (\Gamma )$}, \textit{then}  $T_{\lambda } \delta (\omega
)\in \Lambda (\Omega )$.

Together with Lemma 1.9, this gives

PROPOSITION 1.11.  \textit{If  $Q\in {\rm ex}\, L$},\textit{ then}
$Q:\, K(\Omega )\to \Lambda (\Omega )$.

This claim is used below when proving Proposition 3.2. In
particular, if $L=\{ Q \} $, the $Q: K(\Omega )\to \Lambda (\Omega
)$.

\bigskip

\centerline{ \bf2. Conditions for ordinary and tame dynamics}

\bigskip

In this section, we describe relationships between properties of
the semicascade $(\Omega ,\varphi )$. These properties are mostly
stated in terms of topological characteristics of the associated
semigroups $E(\Omega ,\varphi )$, $\Gamma (\Omega ,\varphi )$,
$G(\Omega ,\varphi )$. The base property is that of $(\Omega
,\varphi )$ being ordinary.

DEFINITION 2.1. A semicascade $(\Omega ,\varphi )$ is said to be
\textit{ordinary} if its enveloping semigroup $E(\Omega ,\varphi )$
is metrizable.

It turns out that the metrizability of the enveloping semigroup
$E(\Omega ,\varphi )$ with the \textit{$\tau $-}topology of the
Cartesian product implies the uniqueness of the representing
measure in Choquet's formula (1.4). In turn, the latter property
guarantees that the semigroups $\Gamma (\Omega ,\varphi )$ and
$E(\Omega ,\varphi )$ essentially coincide and that $\Gamma (\Omega
,\varphi )$ coincides with the set ${\rm ex}\, G(\Omega ,\varphi )$
of extreme points. One more characteristic of the ordinary
semicascade $(\Omega ,\varphi )$ is related to a slight restatement
of the version in [10--12] of the notion of nonchaotic dynamics.

DEFINITION 2.2.  We say that a dynamical system  $(\Omega ,\varphi
)$ is \textit{nonchaotic} if each closed semi-invariant subset
$(\varphi \Theta \subset \Theta )$ of the set $\Theta \subset
\Omega $ contains a trajectory $\, o(\omega ), \, \omega \in \Theta
$, that is Lyapunov stable with respect to the semicascade $(\Theta
,\varphi )$.

Here $o(\omega ) =\{  \varphi ^{n} \omega ,\, \, n\ge 0\, \} $.
Consider yet another type of dynamical systems starting from the
statement due to K\"{o}hler [5].

DEFINITION 2.3.  We say that a semicascade  $(\Omega ,\varphi )$ is
\textit{tame} if, for each continuous function  $x\in C(\Omega )$
and an arbitrary increasing subsequence $\{ n(k)\} $ of positive
integers, one has the relation
       \[\mathop{\inf \, }\limits_{a} \left\| \,
       \sum _{k=1}^{\infty }a_{k} x_{n(k)}  \,
       \right\| _{C(\Omega )}  = \; 0,  \]
where $x_{n} (\omega )=x(\varphi ^{n} \omega )$, $\omega \in
\Omega$, and the infimum is taken over finite real sequences $a=\{
a_{k} \} $ such that $\sum _{k=1}^{\infty}\left|a_{k} \right| \,
=\, 1$.

Both classes of dynamical systems are stable (see [10]) with
respect to passage to closed subsystems. Furthermore, it is well
known [10] that a dynamical system is tame if and only if $E(\Omega
,\varphi )\subset \Pi _{1} $, where $\Pi _{1} $ is the set of first
Baire class transformations of $\Omega $. Let $E_{1} (\Omega
,\varphi )$ be the sequential $\tau $-closure of the set $E_{0}
(\Omega ,\varphi )=\{  \varphi ^{n} ,\, n\ge 0 \} $. Clearly,
$E_{1} (\Omega ,\varphi )\subset E(\Omega ,\varphi )\, \bigcap \,
\Pi _{1} $. Recall that $\Pi _{b}$ is the set of Borel mappings $p:
\Omega \to \Omega $. For brevity, we sometimes write $\Gamma
(\Omega ,\varphi )=\Gamma $, $G(\Omega ,\varphi )=G$.  Just as
above, by $A(\cdot )$ and $K(\cdot )$ we denote the spaces, compact
in the $w^{*} $-topology, of regular Borel probability measures and
Dirac measures, respectively, on some Hausdorff compact set.

Now let us present a table of properties, divided into classes
$D1$--$D6$, of the semicascade $(\Omega ,\varphi )$.

\noindent $D1$:  (a1)  The compact set  $E(\Omega ,\varphi )$ is
metrizable; (b1) the system  $(\Omega ,\varphi )$ is nonchaotic;
(c1) the compact set $G(\Omega ,\varphi )$ is metrizable.

\noindent $D2$: (a2) The compact set $G(\Omega ,\varphi )$ is a
Fr\'echet--Urysohn space; (b2)  ${\rm card}\, G(\Omega
,\varphi)=\mathfrak {c} $.

\noindent $D3$:  (a3)  The semicascade  $(\Omega ,\varphi )$ is a
tame system; (b3)  the compact set  $E(\Omega ,\varphi )$ is a
Fr\'echet--Urysohn space; (c3)  ${\rm card}\, E(\Omega ,\varphi )
\, \le \, \mathfrak {c} $; (d3)  $E(\Omega ,\varphi )=E_{1} (\Omega
,\varphi )$;   (e3)  $\; E(\Omega ,\varphi )\subset \Pi _{b} $.

\noindent $D4$:  (a4) The operators $T\in G(\Omega ,\varphi )$ are
determined by their values on Dirac measures $\mu \in K(\Omega)$.

\noindent $D5$:  (a5)  The convex operator set $G(\Omega ,\varphi
)$ is a Choquet simplex.

\noindent $D6$:  (a6)  The operators  $P\in \Gamma (\Omega ,\varphi
)$ are determined by their values on Dirac measures $\mu \in
K(\Omega)$; (b6)  The canonical epimorphism  $\pi : \Gamma (\Omega
,\varphi )\to E(\Omega ,\varphi )$ is injective; (c6)  ${\rm ex}\,
G=\Gamma $;   (d6)  ${\rm ex}\,  \Gamma =\Gamma $.

Part of the properties of the dynamical system $(\Omega ,\varphi )$
in the classes $D1$, $D3$, and $D6$ were considered earlier in
[10--13]. Although these papers deal with invertible dynamical
systems, the passage to the noninvertible case proves to be purely
technical in situations of interest. Properties (c6) and (d6) were
discussed in [1]. A topological compact set has the
Fr\'echet--Urysohn property if the operation of closure for sets in
it can be defined in terms of sequences. Metrizability implies this
property, but the converse is not true.

Property (a5) means [15, Chap. 9] that the correspondence $\lambda
\stackrel{\gamma }{\longrightarrow}T_{\lambda } $ in the Choquet
representation (1.4) is injective,  i.e., the  cone in ${\rm End}\,
X^{*} $ with base $G(\Omega ,\varphi )$ is a vector structure. If
condition (b6) holds, then the function $\pi : \Gamma (\Omega
,\varphi )\to E(\Omega ,\varphi )$ defined by (1.2) is
simultaneously a homeomorphism and an algebraic isomorphism, which
permits completely identifying the enveloping semigroups $\Gamma
(\Omega ,\varphi )$ and $E(\Omega ,\varphi )$. In this situation,
we write $\Gamma (\Omega ,\varphi )\simeq E(\Omega ,\varphi )$.

THEOREM 2.4.  \textit{The properties of the semicascade  $(\Omega
,\varphi )$ within each of the classes  }D1, D2, D3, and D6\textit{
are pairwise equivalent}.

Thus, each of properties (b1) and (c1) can be taken for an
alternative definition of being ordinary, and each of properties
(b3), (c3), (d3), and (e3) completely characterizes tame systems
$(\Omega ,\varphi )$. We point out that the implications (a1)
$\Leftrightarrow $ (b1) and (a3) $\Leftrightarrow $ (b3)
$\Leftrightarrow $ (c3) were established in [10--13], while the
implications (b6) $\Leftrightarrow $ (c6) $\Leftrightarrow $ (d6)
were obtained in [1]. The equivalence of properties (d3) and (e3)
means that if the enveloping semigroup $E(\Omega ,\varphi )$
contains transformations that do not belong to the first Baire
class, then it also contains transformations that are not Borel.

PROOF.  The class $D1$.  By [10, Secs.~8 and~9],  conditions (a1)
and (b1) are equivalent.  The metrizability of the compact set
$G(\Omega ,\varphi )$ implies that of its closed subset $\Gamma
(\Omega ,\varphi )$ and hence [16] the metrizability of the range
of the continuous mapping $\Gamma (\Omega ,\varphi )\stackrel{\pi
}{\longrightarrow}E(\Omega ,\varphi )$. Thus,  (c1)
$\Rightarrow $  (a1).  Next,  (a1)  $\Rightarrow $  (c3), and
condition (a1) implies the metrizability of the topological space
$A(E)$ of measures. It is well known [10, Theorem 7.1] that (c3)
$\Rightarrow $ (b6), and since, under this condition, $\Gamma
(\Omega ,\varphi )\simeq E(\Omega ,\varphi )$, it follows that
condition (a1) also implies the metrizability of the space
$A(\Gamma )$ of measures. The operator semigroup $G(\Omega ,\varphi
)$ is the range of the continuous surjection $A(\Gamma
)\stackrel{\gamma }{\longrightarrow}G(\Omega ,\varphi )$, and hence
the compact set $G(\Omega ,\varphi )$ is metrizable as well, and
(a1) $\Rightarrow $ (c1).

The class $D2$. As was already noted, the semigroup  $G(\Omega
,\varphi )$ is enveloping for the action  $W\times A\, \,
\stackrel{V}{\longrightarrow}\, A$ on the compact set $A=A(\Omega
)$ of the abelian semigroup $W$ of nonnegative finite numerical
sequences with unit sum and with convolution as multiplication.
Thus, the equivalence of conditions (a2) and (b2) follows from the
results in [10, Sec.~6] for the enveloping semigroup $G=E(A,W)$.

The class $D3$.  The pairwise equivalence of conditions (a3)--(c3)
is a well-known fact [10, Sec.~6], and the implications (b3)
$\Rightarrow $ (d3) $\Rightarrow $ (e3) are trivial.  To establish
the implication (e3) $\Rightarrow $(c3), we slightly generalize the
scheme in [17, Sec.~2.12(ii)] and identify Borel transformations
$p: \Omega \to \Omega $ with their graphs $(\omega, p\, \omega )
\subset \Omega \times \Omega $, $\omega \in   \Omega $. The
transformation $(\omega _{1}, \omega _{2} )\to (\omega _{1} ,p \,
\omega _{2} )$ of the metrizable compact set $\Omega \times \Omega
$ proves to be a Borel transformation, and hence so is the scalar
function $g(\omega _{1} ,\omega _{2} )=\rho (\omega _{2} , p\,
\omega _{1} )$ on $\Omega \times \Omega $ for any of the metrics
$\rho $ compatible with the original topology on $\Omega $. The
graph of $p$ coincides with $g^{-1} (0)$ and hence is a Borel
subset of $\Omega \times \Omega $. It is well known (e.g., see
[16]), that the set of such sets is at most of the cardinality of
continuum. Thus, if $E(\Omega ,\varphi ) \subset \Pi _{b} $, then
${\rm card}\, E(\Omega ,\varphi ) \le \mathfrak {c} $, as desired.

The class $D6$.  Property (a6)  is a useful restatement of (b6).
The equivalence of conditions (b6) and (c6) was established in [1,
Sec.~4]. It was also shown there that one always has ${\rm ex}\,
G={\rm ex}\, \Gamma $. Hence (c6) $\Leftrightarrow $  (d6), and the
proof of the theorem is complete.

Thus, now we can speak of the classes  $\mathcal {D}$1--$\mathcal
{D}$6 of discrete dynamical systems $(\Omega ,\varphi )$ possessing
any of the properties in the respective classes D1--D6. Note that
$\mathcal {D}$1 is the class of ordinary semicascades, and
$\mathcal {D}$3 is the class of tame semicascades.

THEOREM 2.5.  \textit{One has the embeddings}  $\mathcal {D}$1
$\subset \mathcal {D}2$ $\subset \mathcal {D}$3 $\subset \mathcal
{D}$6 \textit{and} $\mathcal {D}$1 $\subset \mathcal {D}4$ $\subset
\mathcal {D}$5 $\subset \mathcal {D}$6.

PROOF.  It suffices to establish the logical chains (a1)
$\Rightarrow $ (b2) $\Rightarrow $ (c3) $\Rightarrow $ (b6)  and
(a1) $\Rightarrow $ (a4) $\Rightarrow $ (a5) $\Rightarrow $ (c6).
It has already been mentioned above that (c3) $\Rightarrow $ (b6)
and shown that (a1) $\Leftrightarrow $ (c1). The implication (c1)
$\Rightarrow $ (b2) is trivial, and since the set $E(\Omega
,\varphi )$ is a continuous image of the compact set $\Gamma
(\Omega ,\varphi )$ and $\Gamma (\Omega ,\varphi )\subset G(\Omega
,\varphi )\; $, it follows that (b2) $\Rightarrow $ (c3),  and
hence (a1) $\Rightarrow $ (b2) $\Rightarrow $ (c3), which proves
the first part of the theorem.

Next, consider the following subsets of the space $C(\Gamma )$ of
scalar functions continuous on the compact set $\Gamma (\Omega
,\varphi )$:
  \[C_{0} (\Gamma )=\{g\in C(\Gamma ): g(P)
  =(x,P\mu ),\; x\in X,\; \mu \in K(\Omega )\} \]
and
  \[C_{1} (\Gamma )=\{g\in C(\Gamma ):g(P)
  =(x,P\mu ),\; x\in X,\; \mu \in A(\Omega )\} .\]
Recall that $X=C(\Omega)$. Properties (a4) and (a5) are equivalent
to saying that the linear spans ${\rm Lin}\, \, C_{0} (\Gamma )$
and ${\rm Lin}\, \, C_{1} (\Gamma )$, respectively, are dense in
the uniform topology on $C(\Gamma )$. Thus,  (a4) $\Rightarrow $
(a5), because $C_{0} (\Gamma )\subset C_{1} (\Gamma )$. Now let us
show that (a1) $\Rightarrow $ (a4). We use the already known
implications (a1) $\Rightarrow $ (c3) $\Rightarrow $ (b6) and (c3)
$\Rightarrow $ (e3). We rewrite formula (1.5) for the operator
$T\in G(\Omega ,\varphi )$ in the form
    \[(T^{*} x,\mu )=\int _{\Gamma }(P^{*} x,\mu )\,
    \lambda (dP) \]
for any $x\in X$ and $\mu \in A(\Omega )$ and note once more that
the function $g(P)=(P^{*} x,\mu )=(x,P\mu )$ is continuous on
$\Gamma (\Omega ,\varphi )$. Since $\Gamma (\Omega ,\varphi )\simeq
E(\Omega ,\varphi )$ in this case and the compact topological space
$E(\Omega ,\varphi )$ is metrizable by property (a1), it follows
that the compact set $\Gamma (\Omega ,\varphi )$ is metrizable as
well, and hence $\lambda _{k} \stackrel{w^{*}
}{\longrightarrow}\lambda $ for some sequence
$\lambda _{k} $ of convex combinations of Dirac measures on $\Gamma
$. Thus,
      \[(T^{*} x,\mu )=\mathop{\lim }\limits_{k\to \infty }
      \, \, \int _{\Gamma }(P^{*} x,\mu )\,
      \lambda _{k} (dP) =\mathop{\lim }\limits_{k\to \infty }
      \, (z_{k},\mu ),  \]
where $z_{k} =\, \int _{\Gamma }P^{*} x\, \lambda _{k} (dP)$,
$z_{k} \in X^{**}$. By using property (b6) and Proposition 1.5, we
find that all operators $P\in \Gamma (\Omega ,\varphi )$ are
regular and, by Lemma 1.2, $P^{*} x\in X_{b} $ whenever $x\in X$.
Thus, $z_{k} \in X_{b} $,  and since one can set $\mu =\delta
(\omega )$, $\omega \in \Omega $, in formula (1.5), we see that
$z_{k} \to T^{*} x$ pointwise on $\Omega $ and $T^{*} x\in X_{b} $
for each $x\in X$. Now the identity $(T^{*} x,\mu )=(x,T\mu )$
implies that all operators $T\in G(\Omega ,\varphi )$ are
determined by their values on the Dirac measures $\delta (\omega )$
and that the semicascade $(\Omega ,\varphi )$ has property (a4);
i.e., (a1) $\Rightarrow $ (a4).

It remains to show that (a5)  $\Rightarrow $  (c6). Assume that
$T\in \Gamma (\Omega ,\varphi )$ but  $T\notin {\rm ex}\, G$; then,
by [18, Theorem 10.1.7], there exists a measure $\lambda \in
A(\Gamma )$ representing the operator $T$ and different from the
Dirac measure $\delta (T)$. This, however, contradicts condition
(a5),  and we have confirmed the logical chain (a1) $\Rightarrow $
(a4) $\Rightarrow $ (a5) $\Rightarrow $ (c6).  The proof of Theorem
2.5 is complete.

The papers [10] and [12] give examples of compact tame but
nonordinary dynamical systems that are however more complicated
than cascades or semicascades. According to [11], the semicascades
generated by homeomorphisms of the segment or the circle are
ordinary. On the other hand [11], an arbitrary semicascade
corresponding to a hyperbolic diffeomorphism of the two-dimensional
torus is not ordinary. Finally, there exists an example [10,
p.~2356] of a minimal distal cascade on the two-dimensional torus
such that this cascade does not belong even to the widest (of the
described above) class $\mathcal {D}$6 and hence is not tame. The
classification of semicascades presented here arose solely in
connection with the problems of ergodic theory considered in the
next section and is of preliminary character. It would be useful to
find out whether all classes $\mathcal {D}$1--$\mathcal {D}$6 of
discrete dynamical systems are pairwise distinct and whether they
form an increasing chain.

Theorem 2.4 shows that a tame semicascade $(\Omega ,\varphi )$
satisfies the relations $E(\Omega ,\varphi )=E_{1} (\Omega ,\varphi
)=E_{b} (\Omega ,\varphi )$. Thus, using also the implication (a3)
$\Rightarrow $ (b6), from Proposition 1.5 we derive the following
important corollary.

COROLLARY 2.6.  \textit{For a tame semicascade $(\Omega ,\varphi
)$, all operators $P\in \Gamma (\Omega ,\varphi )$ are regular}.

The following assertion shows that the enveloping semigroups of
nonordinary dynamical systems have a rather intricate structure.

PROPOSITION  2.7.  \textit{If for a semicascade  $(\Omega , \varphi
)$ there exists a countable everywhere dense set $e\subset \Omega $
such that, for any transformations $p_{1}, p_{2} \in E(\Omega
,\varphi )$, their coincidence on this set implies that $p_{1}
=p_{2} $, then the semigroup $E(\Omega ,\varphi )$ is metrizable}.

PROOF. The space $C(\Omega )$ of scalar continuous functions on a
metrizable compact set $\Omega $ is separable. Let $\{ x_{n} \}$ be
a set of functions countable and norm dense in the closed unit ball
of the space $X=C(\Omega )$, and let $\{ \omega _{m} \}$ be a
sequence of points dense in $\Omega $. For elements $p_{1} ,p_{2}
\in E(\Omega ,\varphi )$, set
  \[d(p_{1} ,p_{2} )=\sum _{n,m=1}^{\infty }\, \frac{1}{2^{n+m} }  \, \,
  \left|\, x_{n} (p_{1} \omega _{m} )-x_{n} (p_{2} \omega _{m}
  )\right|.\]
If $p_{1} \omega _{m} \ne p_{2} \omega _{m} $ for some $m \ge  1$,
then there exists a continuous function $x_{n} $ such that $x_{n}
(p_{1} \omega _{m} )\ne x_{n} (p_{2} \omega _{m} )$, and so
$d(p_{1} ,p_{2} )>0$. Thus, it follows from the assumptions in the
proposition that $p_{1}=p_{2}$ whenever $d(p_{1} ,p_{2} )=0$, and
hence the function $d$ is a metric on the $\tau $-compact set
$E(\Omega ,\varphi )$. The topology generated by this metric is
formally weaker than the $\tau $-topology of pointwise convergence
on $E(\Omega ,\varphi )$. Consequently  [14, Theorem 1.5.8],  these
topologies coincide, and the proof is complete.

We see that in the nonordinary case, for each countable (not
necessarily dense) subset $e\subset \Omega $, the enveloping
semigroup $E(\Omega ,\varphi )$ contains transformations $p_{1} \ne
p_{2} $ coinciding on $e$. On the other hand, Proposition 2.7 has
the following corollary

COROLLARY 2.8.  \textit{If all transformations $p\in E(\Omega
,\varphi )$ are continuous outside some given countable set
$e\subset \Omega $, then the semicascade $(\Omega ,\varphi )$ is
ordinary.}

It is of interest to compare this assertion with the well-known
criterion [19] for the continuity of all elements of the enveloping
semigroup: $E(\Omega ,\varphi )\subset C(\Omega ,\Omega )$ if and
only if the dynamical system  $(\Omega ,\varphi )$ is weakly almost
periodic. The latter means that the set of translates $\{ {\kern
1pt} x(\varphi ^{n} \omega ), \, n\ge 0{\kern 1pt} \} $ of an
arbitrary continuous function $x(\omega )$ is relatively compact in
the weak topology of the space $C(\Omega )$.

\bigskip

\centerline { \bf3.  Ergodic properties of semicascades}

\bigskip

The constructions in the preceding sections are used here to study
ergodic properties of ordinary and tame discrete dynamical systems.
When doing so, we heavily use the results in [1]. Just as before,
on the space $X$ of continuous functions on a metrizable compact
set $\Omega $ we consider the shift operator $U$ corresponding to a
continuous transformation $\varphi : \Omega \to \Omega $. Let
$V=U^{*} $. In the operator algebra ${\rm End}\, X^{*} $, we single
out the semigroups $\Gamma _{0} =\{ V^{n} ,\, n\ge 0 \} $ and
$G_{0} = {\rm co}\, \Gamma _{0} $ as well as the compact semigroups
$G=\overline{\rm co} \, \Gamma _{0}$ and $\Gamma
=\overline{\Gamma}_{0}$, where $\overline{(\,\cdot \,)}$ is the
operation of ${\rm W}^{*} {\rm O}$-closure. The notation of these
semigroups will include the dependence on $\Omega $ and $\varphi $
where necessary. Next, let $\Lambda (\Omega )$ be the set of all
Borel probability $\varphi $-ergodic measures on $\Omega $.

Throughout the following, the convergence of nets in ${\rm End}\,
X^{*} $ is understood in the sense of the ${\rm W}^{*} {\rm
O}$-topology. For an arbitrary semicascade $(\Omega ,\varphi )$, we
say that an operator net $T_{\alpha} \in G_{0} (\Omega ,\varphi )$
is ergodic if $T_{\alpha } (I-V)\; \, \stackrel{{\rm W}^{*} {\rm
O}}{\longrightarrow}\; \, 0$, where $I={\rm id}$. Note that
$T_{\alpha} \; =\; R_{\alpha }^{*} $, where $R_{\alpha } \; \in
{\rm End}\, X$. The sequence
  \[V_{n}  \, = \, \frac{1}{n} \, (I+V+...+V^{n-1} ) \]
of Ces\`aro means is ergodic, because $V_{n} (I-V) = n^{-1}
(I-V^{n} )$ and $\left\| \, V_{n} (I-V)\, \right\| \, \le \,
2n^{-1} $. By [1, Lemma 1.1], the kernel $L$ of the semigroup
$G(\Omega ,\varphi )$ is determined by relation (1.1). Moreover [1,
Lemma 1.2], every element in $L$ is the limit of some operator
ergodic net, and conversely, if an ergodic net $T_{\alpha } $
converges to an operator $Q$, then $Q\in L$. Recall that the kernel
$L$ consists of projections $Q\in {\rm End}\, X^{*} $.

THEOREM 3.1.  \textit{If a semicascade }\textbf{\textit{$(\Omega
,\varphi )$}}\textit{  is ordinary}, \textit{ then every ergodic
operator net $T_{\alpha} $ contains a convergent ergodic
subsequence} $T_{\alpha (k)} $.

PROOF.  Let $Q\in L$. Then $T_{\alpha } \stackrel{{\rm W}^{*} {\rm
O}}{\longrightarrow}\, \, Q$ for some ergodic operator net
$T_{\alpha} \in G_{0} $. Since the semicascade \textbf{$(\Omega
,\varphi )$} is assumed to be ordinary, it follows by Theorem 2.5
that the compact set $G(\Omega ,\varphi )$ has the
Fr\'echet--Urysohn property. Consequently, the net $T_{\alpha} $
contains a convergent sequence $T_{\alpha (k)}$ such that $\alpha
(k+1)>\alpha (k)$ in the sense of the partial order on the indexing
set $\{ \alpha  \} $, and the proof of the theorem is complete.

Thus, in the case of ordinary dynamics, for each projection $Q\in
L$ there exists an ergodic sequence of operators $T_{n} \in G_{0} $
converging to $Q$. We see that actually this result (as well as
Theorem 3.1) holds for the class $\mathcal {D}$2 of compact
discrete dynamical systems, which is apparently wider than the
class $\mathcal {D}$1 of ordinary semicascades.

As a straightforward consequence of the last theorem and
Proposition 1.11, we obtain the following assertion.

PROPOSITION 3.2. \textit{If the semicascade
}\textbf{\textit{$(\Omega ,\varphi )$}}\textit{ is ordinary},
\textit{ then for each projection $Q\in {\rm ex}\, L$ there exists
a convergent ergodic sequence of operators $T_{n} \to Q$ such that,
for each $\omega \in \Omega $, the measures $T_{n} \delta (\omega
)$ converge as $n\to \infty $ in the $w^{*}$-topology to the
ergodic measure }$\mu =Q\delta (\omega )$.

In the conventional terminology [20, Chap. 4], this means that the
asymptotic (with respect to the corresponding operator ergodic
sequence $\{ \, T_{n} \, \} $) distribution of trajectories of the
dynamical system \textbf{$(\Omega ,\varphi )$} is determined by
ergodic measures. If, in the ordinary case, ${\rm card}\, L=1$,
then, for every ergodic operator net $T_{\alpha} \in G_{0} $ and
for each $\omega \in \Omega $, one has the convergence $T_{\alpha
} \delta (\omega )\stackrel{w^{*} }{\longrightarrow}\mu _{\omega }
$, where $\mu _{\omega } \in \Lambda (\Omega )$.

The convergence of the ergodic operator net $T_{\alpha } =R_{\alpha
}^{*} $  in the ${\rm W}^{*} {\rm O}$-topology of the space ${\rm
End}\, X^{*} $ implies the pointwise convergence on $\Omega $ of
the net $R_{\alpha} x$ of functions for an arbitrary continuous
function $x\in X$. In the sequential case, the converse assertion
follows from the Lebesgue dominated convergence theorem, but this
is no longer the case in the general situation. However, it turns
out that, when studying the convergence of ergodic operator nets
$T_{\alpha } $ for ordinary dynamical systems, we can restrict
ourselves to verifying the pointwise convergence of the functions
$R_{\alpha}^{} x$ continuous on $\Omega $ for all $x\in X$.

PROPOSITION 3.3. \textit{If the semicascade}\textbf{\textit{$
(\Omega ,\varphi )$}} \textit{is ordinary, then a sufficient
condition for the ${\rm W}^{*} {\rm O}$-convergence of the ergodic
operator net $T_{\alpha } \in G_{0} $ is given by the pointwise
convergence on } \textbf{\textit{$\Omega $}} \textit{of the
function nets $R_{\alpha } x$ for all  $x\in X$}.

The proof readily follows from the inclusion $\mathcal {D}$1
$\subset \mathcal {D}$4 in Theorem 2.5, so that actually the class
$\mathcal {D}$4 of discrete dynamical systems has this property.

In what follows, by $M(\Omega )$ and $Z(\Omega )$ we denote the
union of all minimal sets and the minimal attraction center,
respectively, of the semicascade $(\Omega ,\varphi )$. Recall that
$Z(\Omega )$ coincides with the closure of the union of supports of
all measures $\mu \in \Lambda (\Omega )$, and $M^{c} (\Omega
)\subset Z(\Omega )$, where $(\cdot )^{c}$ is the operation of
closure in $\Omega $; in the general case, the inclusion may be
proper. As was already noted in Sec.~2, the enveloping semigroup
$E(\Omega ,\varphi )$ for a tame semicascade $(\Omega,\varphi )$
coincides with the transformation class $E_{1} (\Omega ,\varphi )$.
This is even more than sufficient for Theorem 2.4 in [1] to hold,
which implies that the sets $M^{c} (\Omega )$ and $Z(\Omega )$
coincide in this case.

PROPOSITION 3.4.  \textit{If the semicascade $(\Omega ,\varphi )$
is a tame system}, \textit{then $M^{c} (\Omega )=Z(\Omega )$.}

Note that the relation $M^{c} (\Omega )=Z(\Omega )$ does not mean
yet that all ergodic measures $\mu \in \Lambda (\Omega )$ are
concentrated on minimal sets.

For a point $\omega \in \Omega$, by $\overline{o}(\omega )$ we
denote the closure of the trajectory $\{\varphi ^{n} \omega ,\,
n\ge 0 \} $ in the topology of $\Omega $.

Slightly varying the statements in [1], consider the following
properties of an arbitrary semicascade $(\Omega ,\varphi )$.

\noindent (A)  $E(\Omega ,\varphi )=E_{1} (\Omega ,\varphi )$.

\noindent (A$_{1} $) The support of every Borel probability
$\varphi $-ergodic measure is a minimal set.

\noindent (B) Every minimal set $m\subset \Omega $ supports exactly
one Borel probability $\varphi $-invariant measure.

\noindent (C) For each $\omega \in \Omega $, the closure of the
trajectory $\overline{o}(\omega )$ contains a unique minimal set.

\noindent (D) All ergodic operator sequences $T_{n} \in G_{0} $
converge.

\noindent (E) All ergodic operator nets $T_{\alpha } \in G_{0} $
converge.

By [1, Theorem 3.2], for an arbitrary semicascade $(\Omega ,\varphi
)$ one has the implications (E) $\Rightarrow $ (B) + (C), (A$_{1}
$) + (B) + (C) $\Rightarrow $ (D). Furthermore [1, Theorem 1.3],
condition (E) is equivalent to the condition that the kernel $L$ of
the semigroup $G(\Omega ,\varphi )$ is a singleton. Condition (A)
is a restatement of the notion of tame system and is so much the
more true for an ordinary semicascade. It is also well known (see
[8--10]) that, for the case of a tame semicascade, every minimal
set supports exactly one invariant probability measure; i.e.,
condition (B) is satisfied. Moreover, all operators $P\in \Gamma
(\Omega ,\varphi )$ are regular in this case  (see Corollary 2.6).
Since $P\mu =\mu $ for ergodic measures $\mu \in \Lambda (\Omega
)$, it follows that $\mu (p^{-1} e)=\mu (e)$ for all
transformations $\, p \in E(\Omega ,\varphi )$ and Borel sets $e
\subset \Omega $. As the proof of [1, Theorem 3.2] shows, this
suffices for deriving the implication (A) + (C) $\Rightarrow
$(A$_{1} $). Thus, in the tame situation we have (A$_{1} $) + (C)
$\Rightarrow $ (D) and hence eventually (C) $\Rightarrow $(A$_{1}
$) and (C) $\Rightarrow $ (D). Since one always has (E)
$\Rightarrow $ (B) + (C) and Theorem 3.1 implies that (D)
$\Leftrightarrow $ (E) in the case of ordinary dynamics, it follows
that the equivalences (C) $\Leftrightarrow $ (D) $\Leftrightarrow $
(E) hold in this class of discrete dynamical systems. Let us
summarize the preceding in the form of two theorems.

THEOREM 3.5. \textit{If the semicascade} $(\Omega ,\varphi
)$\textit{ is ordinary}, \textit{ then all ergodic operator nets
$T_{\alpha} $ converge if and only if, for each $\omega \in \Omega
$, the closure of the trajectory }$\overline{o}(\omega )$
\textit{contains a unique minimal set. The same is true with nets
replaced by sequences}.

THEOREM 3.6. \textit{If, for a tame semicascade }$(\Omega ,\varphi
)$\textit{, the closure of each trajectory }$\overline{o}(\omega
)$\textit{ contains a unique minimal set, then all ergodic operator
sequences $T_{n} $ converge, and the support of each $\varphi
$-ergodic measure is a minimal set.}

We say that a pair $\omega _{1} ,\, \omega _{2} \in \Omega $ of
points is proximal and write $(\omega _{1} ,\, \omega _{2} )\in
{\rm P} $ if $\mathop{\inf \,}\limits_{n\ge 0} \rho (\varphi ^{n}
\omega _{1} ,\, \varphi ^{n} \omega _{2} )=0$ for some given metric
$\rho $ compatible with the given topology on $\Omega $. By [4,
Proposition 5.16], the proximality relation  ${\rm P} ={\rm P}
(\Omega ,\varphi )$ is transitive (i.e., is an equivalence
relation) if and only if the kernel of the enveloping semigroup
$E(\Omega ,\varphi )$ is a minimal left ideal. It follows that the
transitivity of the relation ${\rm P} (\Omega ,\varphi )$ is
independent of the choice of an admissible metric on $\Omega $ .

LEMMA 3.7.  \textit{If the proximality relation ${\rm P} (\Omega
,\varphi )$ is transitive},\textit{ then for each point $\omega \in
\Omega $ there exists exactly one minimal set $m\subset
\overline{o}(\omega )$}.

PROOF. Assume that the proximality relation is transitive and there
exists a point $\omega \in \Omega $ such that the closure
$\overline{o}(\omega )$ contains two distinct minimal sets  $m_{1}
$ and $m_{2} $. Then [1, Lemma 3.2] $I_{1} \omega =m_{1} $ and
$I_{2} \omega =m_{2} $ for some minimal left ideals $I_{1} $ and
$I_{2} $ in $E(\Omega ,\varphi )$. By the general theory of
enveloping semigroups [4, Chap. 3; 5], there exist idempotents
$p_{1} \in I_{1} $ and $p_{2} \in I_{2} $ such that $(\omega ,p_{1}
\omega )\in {\rm P} $ and $(\omega ,p_{2} \omega )\in {\rm P} $. At
the same time, however, $(p_{1} \omega ,p_{2} \omega )\notin {\rm
P} $, because the closed sets $m_{1} $ and $m_{2} $ are disjoint.
This contradiction proves the lemma.

By combining Theorems 3.5 and 3.6 with Lemma 3.7, we obtain the
following corollary.

COROLLARY 3.8. \textit{If, for an ordinary semicascade $(\Omega
,\varphi )$, the proximality relation ${\rm P} (\Omega ,\varphi )$
is transitive, then all ergodic operator nets $T_{\alpha} $
converge}. \textit{In the case of tame dynamics, the transitivity
of the relation ${\rm P} (\Omega ,\varphi )$ implies the
convergence of all ergodic operator sequences $T_{n} $}.

Since the convergence of Ces\`aro means is an important problem in
ergodic theory, we single out the following result, which is a
consequence of Theorem 3.1 and Corollary 3.8.

PROPOSITION 3.9. \textit{For an ordinary semicascade $(\Omega
,\varphi )$, the sequence of Ces\`aro means $V_{n} $ always
contains a convergent subsequence $V_{n(k)} $}.\textit{ If, for a
tame semicascade $(\Omega ,\varphi )$, the proximality relation
${\rm P} (\Omega ,\varphi )$ is transitive}, \textit{ then the
operator means $V_{n} $ converge as $n\to \infty $}.\textit{
}

Note that the claims about ordinary semicascades in Propositions
3.2 and 3.9 as well as in Theorem 3.5 and Corollary 3.8 actually
remain valid in the formally broader class $\mathcal {D}$2 of
compact discrete dynamical systems.

\quad

\centerline { \bf References}

\quad

\noindent [1] A. V. Romanov, ``Weak$^{*}$ convergence of operator
means''. \textit {Izv. Math.}, \textbf{75}:6 (2011),

 1165-1183.

\noindent [2] A. V. Romanov, ``Ordinary semicascades and their
ergodic properties''. \textit{Funct. Anal.}

\textit{Appl.},  \textbf{47}:2 (2013), 160-163.

\noindent   [3]  A. V. Romanov, ``Enveloping semigtoups and ergodic properties
of compact

semicascades''. \textit{International Conference  on Differential Equations and
Dynamical}

\textit{Systems} (\textit {Suzdal} \textit {2012}), \textit {Abstracts}, MIAN,
Moscow, 2012, p.151-152.

\noindent  [4]  R. Ellis,  \textit{Lectures on Topological
Dynamics.}  Benjamin, New York, 1969.

\noindent  [5]  A. K\"{o}hler,\textit{  }``Enveloping semigroups
for flows''.  \textit{Proc. Roy. Irish. Acad}.,
 \textit{Sect. A},

 \textbf{95}:2 (1995), 179-191.

\noindent  [6] N. Kryloff,  N. Bogoliouboff,  ``La theorie
generale de la mesure dans son application a

l'etude des systemes dynamiques de la mecanique non lineaire''.
\textit{Ann. of Math}.,

\textbf{38}:1 (1937),  65-113.

\noindent   [7] S. P. Lloyd   ``On the ergodic theorem of Sine''.
\textit{Proc. Amer. Math. Soc.}, \textbf {56} (1976),

121-126.

\noindent  [8]  W. Huang, ``Tame systems and scrambled pairs under
an abelian group action''.

 \textit{Ergodic Theory Dynam. Systems}, \textbf{26}:5 (2006), 1549-1567.

\noindent  [9]  D. Kerr, H. Li,  ``Independence in topological and
$C^{*} $-dynamics''.
        \textit{Math. Ann}.,  \textbf{338}:4

        (2007), 869-926.

\noindent [10]  E. Glasner, ``Enveloping semigroups in topological
dynamics''.  \textit{Topology Appl.},

\textbf{154}:11 (2007),  2344-2363.

\noindent [11]  E. Glasner, M. Megrelishvili,  ``Hereditarily
non-sensitive dynamical systems and

linear representations''. \textit{Colloq. Math}., \textbf{104}:2
(2006),  223-283.

\noindent [12] E. Glasner, M. Megrelishvili, V. V. Uspenskij, ``On
metrizable enveloping

semigroups''. \textit{Israel J. Math}.,  \textbf{164}:1 (2008), 317-332.
\textit{}

\noindent [13] E. Glasner, ``On tame dynamical systems''.
\textit{Colloq. Math}., \textbf{105}:2 (2006), 283-295.

\noindent

\noindent [14]  N. Dunford and J. T. Schwartz,  \textit {Linear
Operators. General Theory. Part I.}

Interscience Publishers, New York - London, 1958.

\noindent [15]  R. R. Phelps, \textit {Lectures on Choquet's
Theorem.} Van Nostrand, Princeton - Toronto

 - New York - London, 1966.

\noindent [16] N. Bourbaki,  \textit{Elements of Mathematics.
General Topology. Parts I and II.} Hermann,

 Paris, 1966.

\noindent [17]  V. I. Bogachev,  \textit{Measure Theory, Vol.I.}
 Springer, Berlin -- Heidelberg -- New York,

2007.

\noindent [18]  R. E. Edwards,  \textit{Functional Analysis.
Theory and Applications}.  Holt, Rinehart and

Winston, New York, 1965.

\noindent [19]  R. Ellis, M. Nerurkar, ``Weakly almost periodic
flows''. \textit{Trans. Amer. Math. Soc.},

        \textbf {313}:1 (1989), 103-119.

\noindent [20]   A. Katok and B. Hasselblatt, \textit
{Introduction to the Modern Theory of Dynamical}

\textit {Systems.} Encyclopedia Math. Appl., vol. \textbf{54},
Cambridge Univ. Press, Cambridge,

1995.

\end{document}